\def\Re{{\rm Re}\,}
\def\halmos{\hbox{\vrule height0.15cm width0.01cm\vbox{\hrule height
 0.01cm width0.2cm \vskip0.15cm \hrule height 0.01cm width0.2cm}\vrule
 height0.15cm width 0.01cm}}
\font\eightrm=cmr8  
\font\eighttt=cmtt8
\magnification=\magstephalf

\parindent=0pt
\overfullrule=0in
\bf
\noindent
A HIGH-SCHOOL ALGEBRA\footnote{$^1$}
{\eightrm 
and high-school (purely formal) calculus.
}
, WALLET-SIZED PROOF, OF THE BIEBERBACH 
\hbox{ \hskip96pt \relax CONJECTURE [After L. Weinstein]}
\medskip
\noindent
\it
\hskip100pt \relax Shalosh B. Ekhad${}^2$ and Doron Zeilberger\footnote{$^2$}
{\eightrm 
Department of Mathematics, Temple University,
Philadelphia, PA 19122, USA. 
\break
{\eighttt ekhad@math.temple.edu ,}
{\eighttt zeilberg@math.temple.edu .}
Supported in part by the NSF. We would like to thank Richard Askey, 
Ira Gessel, Jeff Lagarias and Herb Wilf for comments, that improved readability.
}
\bigskip
\it
\qquad \qquad \qquad \qquad \qquad
\qquad\qquad Dedicated to Leonard Carlitz\footnote{$^3$}
{\eightrm
L. Carlitz was, for many years, editor of the Duke Journal, until
he was relieved from his duties by the proponents of so-called
``modern math'', who despised formal, Carlitz-style, mathematics.
This paper makes explicit the power of formal math, that is (too) subtly
hidden in Weinstein's original presentation.
}
, master of formal mathematics.
\bigskip
\rm
Weinstein's[2]  brilliant short proof of de Branges'[1] theorem can be made
yet much shorter(modulo routine calculations), completely elementary
(modulo L\"owner theory),  self contained(no need for the esoteric
Legendre polynomials' addition theorem), and motivated(ditto),
as follows. Replace the text between p. 62, line 7 and p. 63, line 7, by Fact 1
below, and
the text between the last line of p.63 and p.64, line 7, by Fact 2 below.
\par
\smallskip
\noindent
{\bf FACT 1:} Let $f_t (z) = e^t z \exp ( \sum_{k=1}^{\infty} c_k (t) z^k )$ 
where $c_k (t)$ are {\it formal} functions
of $t$. Let $z$ and $w$ be related by $z/(1-z)^2 = e^t w/(1-w)^2$. The 
following formal identity holds. (For any formal Laurent series $f(z)$,
$CT_z f(z)$ denotes the {\it Constant Term} of $f(z)$.)
\smallskip
$$
(1+w) { {d} \over {dt}} \{
\sum_{k=1}^{\infty} ( 4/k - k  c_k (t) \overline{c_k (t)} ) w^k \} =
$$
$$
(1-w) \sum_{k=1}^{\infty} \Re CT_z
\left\{
{
{ {\partial f_t (z)} \over {\partial t} }
\over
{ z {\partial f_t (z)} \over {\partial z} } 
}
\cdot
\left ( 2(1+  \sum_{r=1}^{k} r c_r (t) z^r ) - k c_k (t) z^k \right ) \cdot
\left ( 2(1+   \sum_{r=1}^{k} r \overline{c_r (t}) z^{-r} ) - 
k \overline{c_k (t)} z^{-k} \right ) 
\right\} w^k
$$
\medskip
\noindent
{\bf Proof:} Routine. (Even for a human.) \halmos
\bigskip
\noindent
{\bf FACT 2:} The coefficients $A_{k,n} (c)$ in the
formal power series (Laurent in $w$) expansion
$(1- z(2c+ (1-c)(w+ 1/w))+z^2 )^{-1}= \sum_{n=0}^{\infty}\sum_{k=0}^{n}
A_{k,n}(c) (w^k + w^{-k}) z^n$ are non-negative for $0 \leq c \leq 1$.
\medskip
\noindent
{\bf Proof:} This follows immediately from the stronger fact that the 
coefficients $B_{k,n} (c)$, defined by the expansion
$(1- z(2c+ (1-c)(w+ 1/w))+z^2 )^{-1/2}=\sum_{n=0}^{\infty}\sum_{k=0}^{n}
B_{k,n}(c) (w^k + w^{-k}) z^n$ can be expressed as $L_{k,n}(c)^2$,
for some double sequence $L_{k,n}$, $0 \leq k \leq n$, that
is {\it real} for  $c$ in $[0,1]$. 
\smallskip
First use Maple to output $B_{k,n}(c)$,
for $0 \leq k \leq n \leq 20$, factor them, and observe that for this
range they are expressible as $L_{k,n}(c)^2$. Using the
{\it gfun} Maple package\footnote{$^4$}
{\eightrm Developed by Salvy and Zimmerman, 
Available by anonymous ftp from {\eighttt
ftp.inria.fr } in directory lang/maple/INRIA/gfun $\quad$ .}
, Maple conjectures
a certain second-order linear recurrence, in $n$, 
satisfied by the $L_{k,n}(c)$,
and using {\it gfun} once again, it computes the linear recurrence
satisfied by the squares of the terms of the solution-sequence of the previous 
recurrence. That new recurrence
turns out to be identical with the third-order recurrence,
in $n$, that the  WZ method\footnote{$^5$}
{\eightrm Wilf and Zeilberger, Invent. Math. 108(1992), 575-633. The
program, and the input file for this problem, are available via
anonymous ftp to {\eighttt math.temple.edu}, 
in directory {\eighttt pub/zeilberger/programs}.}
outputs for $B_{k,n}$. It follows that indeed $B_{k,n}=L_{k,n}^2$,
where {\it now} $L_{k,n}(c)$ denotes the solution of the above-mentioned
second-order recurrence, with the obvious initial values at $n=k,k+1$. \halmos
\medskip
\noindent
{\bf References}
\smallskip
\noindent
{\bf 1.} L. de Branges, {\it A proof of the Bieberbach conjecture}, Acta Math.
{\bf 154}(1985), 137-152.
\par
\noindent
{\bf 2.} L. Weinstein, {\it  The Bieberbach Conjecture}, Inter. Math. Res. Notices
(of the Duke J.) {\bf 3}(1991, No. 5), 61-64. 
\medskip
Revised Version: Nov. 1993. Original Version: Nov. 1992.
\bye